\documentclass{article}
\usepackage{spconf,amsmath,graphicx}

\makeatletter
\def\name#1{\gdef\@name{{\em #1}}} 
\makeatother

\usepackage{xcolor}

\usepackage{graphicx} 
\usepackage{amsthm}
\usepackage[utf8]{inputenc}
\usepackage[table,xcdraw]{xcolor} 
\usepackage{float}
\usepackage{cite}
\usepackage{graphicx}
\usepackage{textcomp}
\usepackage{xcolor}
\usepackage{tikz}
\usepackage{pgfplots}
\usepackage[export]{adjustbox}
\usetikzlibrary{spy}
\usepackage{graphicx}
\usepackage{tabularx}
\usepackage{xcolor}
\usepackage{amssymb,amsmath,amsfonts}
\usepackage{bm}
\usepackage{tikz}
\usepackage{multirow}
\usepackage{booktabs}
\usepackage{blindtext}
\usepackage{multicol}
\usepackage{tikzscale}
\usepackage{numprint}
\npdecimalsign{.}
\nprounddigits{3} 
\usepackage{caption}
\usepackage{subcaption}
\usepackage{enumitem}
\usepackage{float}
\usepackage{pict2e}
\usepackage{siunitx}
\usepackage[english]{babel}
\usepackage{authblk}
\usepackage{orcidlink}
\usepackage{lipsum}
\usepackage{wrapfig}
\usepackage{array}  
\renewcommand{\arraystretch}{2.5}
\usepackage[most]{tcolorbox}
\definecolor{darkspringgreen}{rgb}{0.09, 0.45, 0.27}

\usepackage{mathtools}
\usepackage[percent]{overpic}
\usepackage{algorithm,algcompatible,amsmath}

\DeclareMathOperator*{\argmin}{arg\,min}
\algnewcommand\INPUT{\item[\textbf{Input:}]}%
\algnewcommand\PARAMETER{\item[\textbf{Parameters:}]}%
\algnewcommand\OUTPUT{\item[\textbf{Output:}]}%

\usepackage{algorithm}
\usepackage{algpseudocode}
\algnewcommand{\Inputs}[1]{%
  \State \textbf{Inputs:}
  \Statex \hspace*{\algorithmicindent}\parbox[t]{.8\linewidth}{\raggedright #1}
}
\algnewcommand{\Initialize}[1]{%
  \State \textbf{Initialize:}
  \Statex \hspace*{\algorithmicindent}\parbox[t]{.8\linewidth}{\raggedright #1}
}

\newcommand{\bx}{\bm{x}}
\newcommand{\bh}{\bm{h}}
\newcommand{\bz}{\bm{z}}

\newcommand{\bA}{\bm{A}}
\newcommand{\bb}{\bm{b}}

\title{Split, Skip and Play: Variance-Reduced ProxSkip for Tomography Reconstruction is Extremely Fast}

\name{Evangelos Papoutsellis$^{1}$\qquad Zeljko Kereta$^{2}$\qquad Kostas Papafitsoros$^{3}$}
\address{$^{1}$Finden Ltd\qquad $^{2}$University College London\qquad $^{3}$Queen Mary University of London}

\begin{document}

\maketitle

\begin{abstract}
Many modern iterative solvers for large-scale tomographic reconstruction incur two major computational costs per iteration: expensive forward/adjoint projections to update the data fidelity term and costly proximal computations for the regulariser, often done via inner iterations. This paper studies for the first time the application of methods that couple randomised skipping of the proximal with variance-reduced subset-based optimisation of data-fit term, to simultaneously reduce both costs in challenging tomographic reconstruction tasks. We provide a series of experiments using both synthetic and real data, demonstrating striking speed-ups of the order $5\times$--$20\times$ compared to the non-skipped counterparts which have been so far the standard approach for efficiently solving these problems. Our work lays the groundwork for broader adoption of these methods in inverse problems.

\end{abstract}

\section{Introduction}
\label{sec:intro}

Tomographic  modalities like X-ray Computed Tomography (CT), and Positron Emission Tomography (PET)  reconstruct an unknown image or volume $\bx^{\dagger}\in\mathbb{R}^m$ from indirect and noisy measurements $\bb\in\mathbb{R}^{m}$ produced by a (typically linear but ill-posed) forward operator $\bA:\mathbb{R}^n\rightarrow\mathbb{R}^m$.  
a standard approach is to use variational regularisation and solve
\begin{equation}
  \argmin_{\bx\in\mathbb{R}^n} \mathcal{D}(\bA\bx,\bb) +\alpha\mathcal{R}(\bx).  
  \label{eq:general_opt_problem}
\end{equation}
The data fidelity term $\mathcal{D}$ penalises the mismatch between the predicted measurements $\bA\bm{x}$ and the acquired data $\bb$ and  is determined by the noise statistics, e.g.\ squared $\ell_2$   
under Gaussian noise.  The regulariser $\mathcal{R}$ encodes prior information and promotes solution desirable properties, such as smoothness, sparsity,  or low-rank structure.  Commonly used choices for $\mathcal{R}$   include model based priors such as Total Variation (TV), and its extensions  \cite{Benning_Burger_2018}. Alternatively, one can use learned implicit priors, e.g.\ plug-and-play (PnP) approaches \cite{PnP2013}, where classical proximal regularisation is replaced  by a denoiser, such as BM3D \cite{Dabov2007} or a pretrained neural network.

Problem \eqref{eq:general_opt_problem} falls under the general framework of minimising composite, and potentially non-smooth, objectives
\begin{align}
\min_{\bx\in\mathbb{R}^n} f(\bx) + g(\bx),
\label{eq:fbs_objective}
\end{align}
where $f$ and $g$ correspond to $\mathcal{D}$ and $\mathcal{R}$, respectively.
Proximal Gradient Descent, also known as the Iterative Shrinkage-Thresholding Algorithm (ISTA), and its accelerated variant FISTA \cite{Combettes2005, BeckTeboulle}, are widely used for convex $f$ with $L$-Lipschitz gradient,  and proper, convex $g$,  see Algorithms\ref{alg:all_algos}.

Each iteration of ISTA and FISTA requires evaluating both the gradient of $f$, and the proximal operator of $g$. The latter is defined by
\begin{equation}
\mathrm{prox}_{\tau g}(\bx):=\argmin_{\bz\in\mathbb{R}^n}\bigg\{\frac{1}{2}\|\bz-\bx\|_{2}^{2} + \tau g(\bz)\bigg\}, \quad \tau>0 
\label{eq:proximal_definition}
\end{equation}
and has to be practically tractable.
While these  algorithms are simple and flexible, their practical efficiency is often limited by two dominant per-iteration costs. The first is the gradient evaluation of the data fidelity term which requires repeated application of the forward operator and its adjoint and is costly, especially for large-scale 3D tomography. The second is the proximal step, which can also be expensive, e.g., for TV-type proximals computed via iterative solvers. Moreover, in PnP methods the proximal step is interpreted as a denoising operation which can also be costly e.g., for BM3D. 

\subsection{ProxSkip}

To reduce the computational cost of the proximal operator, a natural strategy is to avoid computing it at every iteration. This idea is formalised by the \emph{ProxSkip} algorithm \cite{mishchenko2022proxskip}, which applies the proximal operator only at randomly selected iterations with probability $p$, while employing a control variate to stabilise the resulting skipping. If $p=1$, the ProxSkip algorithm coincides with ISTA, see Algorithm \ref{alg:all_algos}. Convergence in expectation is established for any $p\in(0,1]$, 
in both the strongly convex and merely convex regimes (ergodic convergence) for the function $f$ \cite{mishchenko2022proxskip, Condat2022}. In practice, ProxSkip is most beneficial when the skipping probability is relatively low. making the proximal evaluations sufficiently infrequent. However, if $p$ is chosen too small, proximal updates are rare and progress can be slow. 
Convergence rate can be optimised by selecting the probability parameter as  $p=\sqrt{\frac{\mu}{L}}$. This, however, relies on strong convexity, a rather restrictive assumption in imaging inverse problems. Moreover,  this choice is not necessarily optimal in terms of computational time to reach a target accuracy.

ProxSkip was originally proposed in the context of federated/distributed optimisation. It has only recently been explored for imaging inverse problems \cite{Papoutsellis2025}, where it was shown to substantially reduce the computational burden and run-times of expensive proximal evaluations, while achieving the same image quality. Nevertheless, the cost of forward/adjoint evaluations in tomography remains a significant burden.

\subsection{Data Spliting and Variance Reduced Estimators}

\begin{figure}
    \centering
    \includegraphics[height=4.3cm]{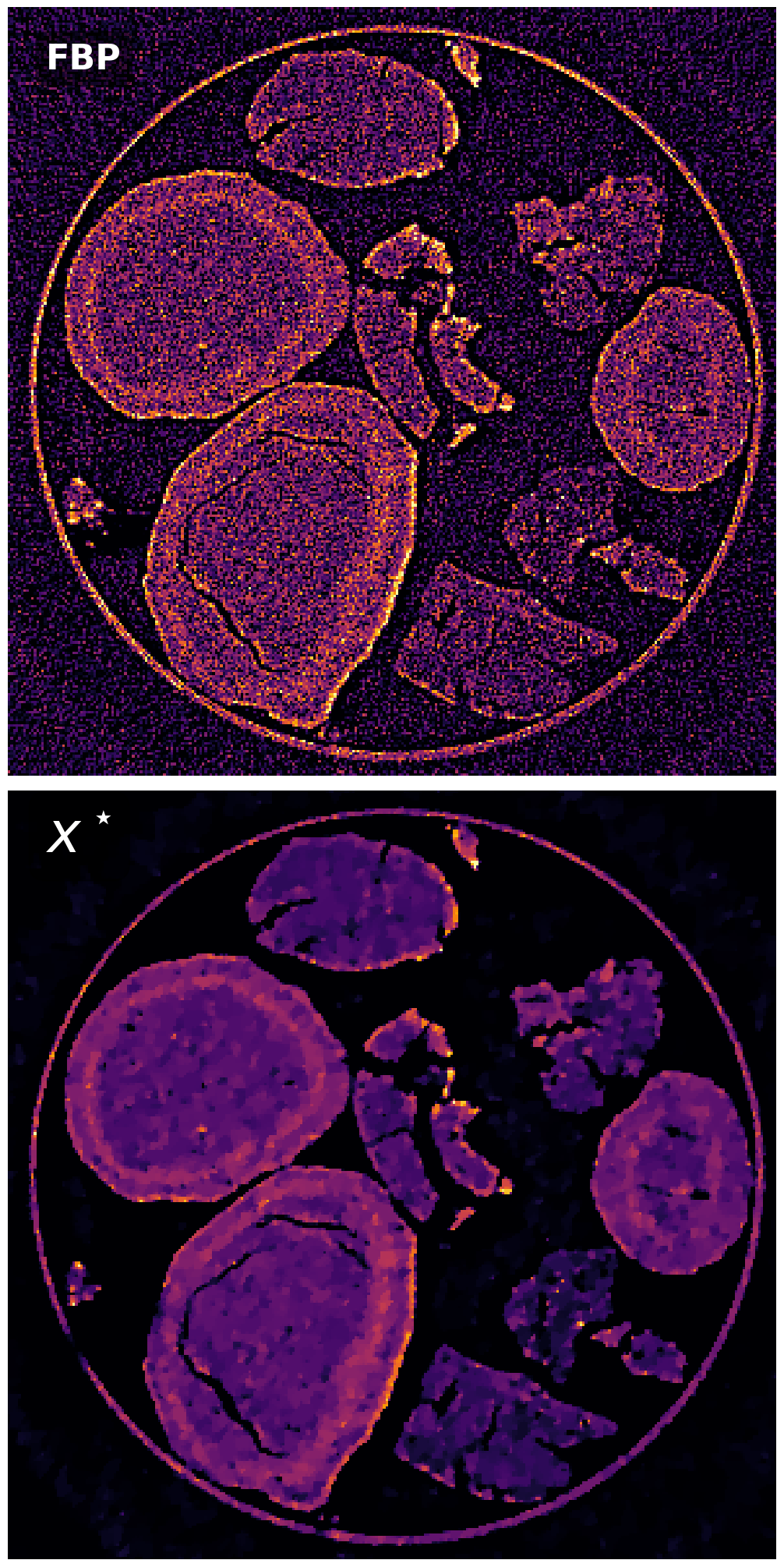}
    \includegraphics[height=4.3cm]{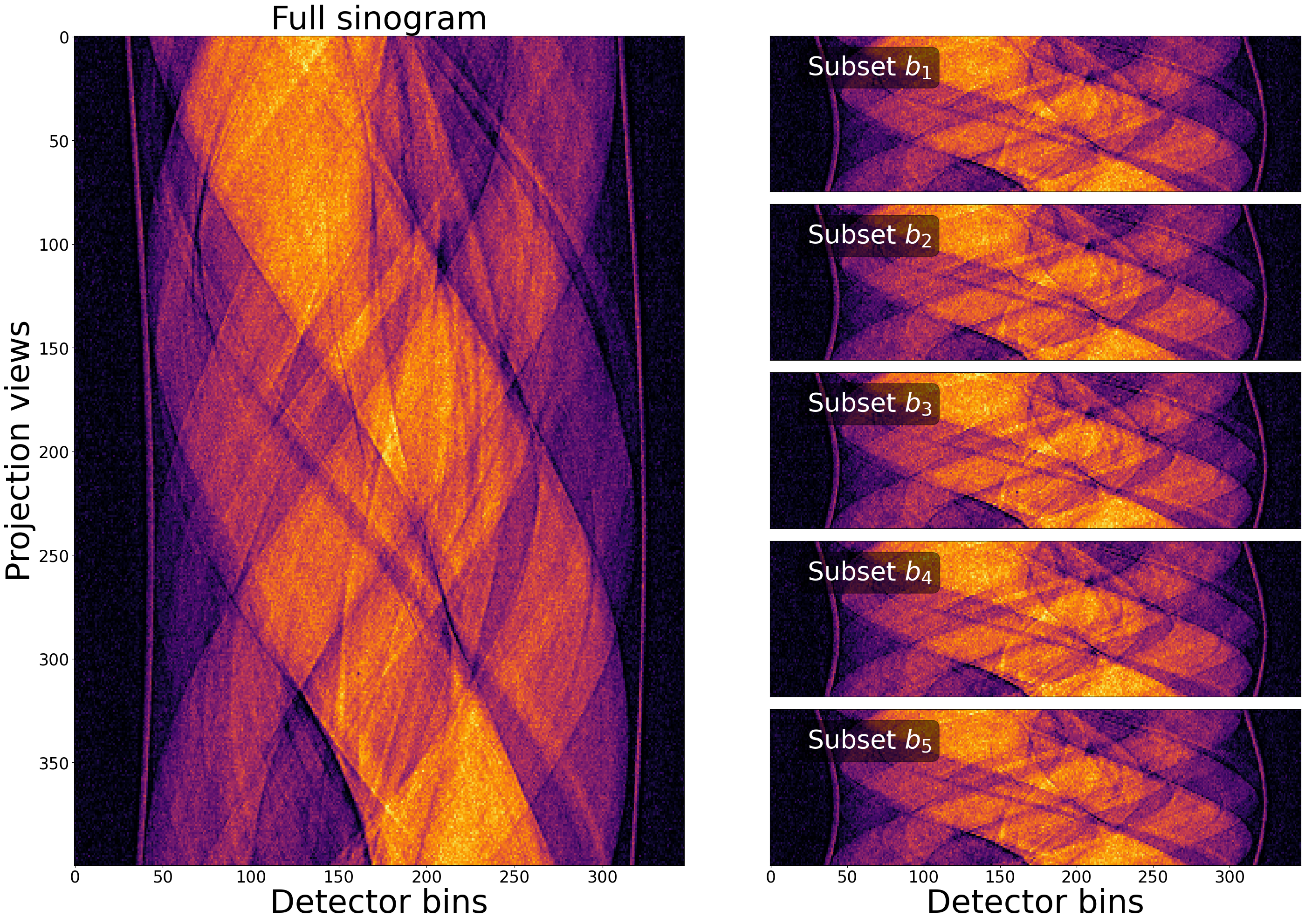}
    \caption{Left: FBP reconstruction  of the post–partial oxidation of methane reaction \cite{Matras2021} and high-accuracy reference (TV-regularised) solution, computed with $2\times 10^{5}$ PDHG iterations using diagonal preconditioning \cite{ChambollePock_2011_Diagonal}. Right: Example of data splitting ($N=5$) using staggered indexing for the projection angles: $b_k=\{k,\,k+5,\,k+10,\,\ldots\}$ for $k=0,\ldots,4$.}
    \label{fig:chemical_dataset_data_splitting}
\end{figure}

In cases where the forward model can be decomposed into blocks, e.g. as subsets of projection angles in CT/PET \cite{Papoutsellis2021}, temporal frames in dynamic imaging \cite{Protopapa2025}, or groups of measurements in computational photography \cite{Sun2019},   the data fidelity can be decomposed across measurements.
Here one can express the data fidelity term in \eqref{eq:general_opt_problem} in a finite-sum form as
$
\mathcal{D}(\bA\bx,\bb) = \sum_{i=1}^{N} D(\bA_i \bx, \bb_i)
$
where $\{(\bA_i,\bb_i)\}_{i=1}^{N}$ represent subsets of the forward operator and data. In tomography, each $\bb_i$ corresponds to a subset of projection angles, see Figure \ref{fig:chemical_dataset_data_splitting}. 

Evaluating $\nabla D(\bA_i\bx,\bb_i)=\nabla f_{i}(\bx)$ typically only requires applying $\bA_i$ and $\bA_i^\top$, i.e.\ forward and backward projections for this subset. 

For example, when applied to CT reconstruction, SGD uses $G_k(\bx)=N \bA_i^\top (\bA_i\bx-\bb_{i})$ as an estimator of the full gradient $\bA^\top(\bA\bx-\bb)$.

Variance introduced around the true gradient may prevent convergence to high accuracy solutions under constant step sizes, causing oscillations or stagnation in ill-conditioned inverse problems. These issues can be addressed with \emph{Variance-Reduction} (VR) techniques, which combine cheap stochastic gradient evaluations while controling the variance. 
 
For example, in Stochastic Average Gradient Amelior\'e (SAGA) \cite{Defazio2014} this is achieved by storing a table of subset gradients
$\{v_k^{\,i}\}_{i=1}^{N}$ and defining $\bar v_k = \sum_{i=1}^{N} v_k^{\,i}$. By sampling $i_k$ the SAGA estimator for the gradient at iteration $k$ is 
\begin{equation}\label{saga_estimator}
\bm{G}_k(\bx) :=
N(\nabla f_{i_k}(\bx)-v_k^{\,i_k})+\bar v_k,
\end{equation}
and the
table is updated only at index $j=i_k$ with $\nabla f_j(\bx)$. 
Stochastic Variance Reduced Gradient (SVRG) \cite{Johnson2013, Xiao2014} instead uses a snapshot point $\tilde{\bx}$ and the full gradient $\nabla f(\tilde{\bx}) := \sum_{i=1}^{N}\nabla f_i(\tilde{\bx})$, which are stored and updated every $N$ iterations. The gradient estimator is computed as 
$$
\bm{G}_k(\bx):=
N(\nabla f_{i_k}(\bx)-\nabla f_{i_k}(\tilde{\bx}))+\nabla f(\tilde{\bx}),
$$
Loopless SVRG \cite{LSVRG} replaces the deterministic schedule for updating the full gradient with a probabilistic one; in each iteration full gradient is updated with probability $1/N$, cf.\ \cite{Richtarik2020}.

\begin{algorithm}[H]
\begin{algorithmic}[1]
\State {\bf Parameters:} $\gamma > 0$, probability $p>0$, data subsets $N$
\State {\bf Initialize:} $\bx_0, \bh_{0} \in \mathbb{R}^n$ 
\For{$k=0,1,\dotsc,K-1$}
\State Compute $\bm{G}_{k}$ (\,unbiased estimator of $\nabla f(\bx_k)\,)$
\State $\hat \bx_{k+1} = \bx_k - \gamma (\bm{G}_{k} (\bx_k) - {\bh_k})$ \hfill 
\State Sample $\theta_{k}\sim\mathrm{Bernoulli}(p), \,\theta_{k}\in\{0, 1\}$ 
\If{$\theta_{k}=1$} 
\State  $\bx_{k+1} = \mathrm{prox}_{\frac{\gamma}{p}g}\bigg(\hat \bx_{k+1} - \frac{\gamma}{p}{\bh_k} \bigg)$ 
\Else
$\;\;\bx_{k+1} = \hat \bx_{k+1}$ \hfill 
\EndIf
\State ${ \bh_{k+1}} = {\bh_k} + \frac{p}{\gamma}(\bx_{k+1} - \hat \bx_{k+1})$ 
\EndFor  
\end{algorithmic}

\begingroup
\begin{center}
\small
\setlength{\tabcolsep}{3pt}
\renewcommand{\arraystretch}{0.8}
\resizebox{0.98\linewidth}{!}{
\begin{tabular}{llc}
\toprule
$p=1$     & $0<p<1$ & $\bm{G}_{k}$   \\
\midrule
ISTA          & ProxSkip   & $\nabla f(\bx_k)$         \\[0.4em]   
ProxSGD      & ProxSGDSkip & $N\nabla f_{i_k}(\bx)$\\ [0.3em] 
ProxSAGA     & ProxSAGASkip & $N(\nabla f_{i_k}(\bx)-v_k^{\,i_k})+\bar v_k$\\[0.4em] 
ProxSVRG   & ProxSVRGSkip & $N(\nabla f_{i_k}(\bx)-\nabla f_{i_k}(\tilde{\bx}))+\nabla f(\tilde{\bx})$\\[0.5em] 
ProxLSVRG   & ProxLSVRGSkip & as above, updated with $p=1/N$ \\[0.4em]
\midrule
\multicolumn{3}{l}{FISTA is ISTA  with an acceleration step \cite{BeckTeboulle}}\\
\bottomrule
\end{tabular}
}
\end{center}
\vspace{-0.5em}
\endgroup

\caption{Family of all considered algorithms}
\label{alg:all_algos}
\end{algorithm}

Combining proximal skipping with stochastic VR gradients yields a family of \emph{ProxSkip-VR} methods (ProxSAGASkip, ProxSVRGSkip,
ProxLSVRGSkip), all outlined in Algorithms \ref{alg:all_algos} that reduce both dominant per-iteration costs in terms of the forward model and the proximal operator.

 Hence, ProxSkip-VR provides two complementary computational savings, \emph{data splitting} for cheaper gradient steps and \emph{proximal skipping} for cheaper regularisation steps, while maintaining stable convergence behaviour comparable to full-gradient proximal methods.
ProxSkip-VR algorithms were only recently introduced for the problem of federated learning \cite{vrproxskip2022}.

\subsection{Our contribution}
\label{sec:contribution}

 Here, we show for the first time the feasibility of ProxSkip-VR algorithms in challenging large scale inverse problems. We use the proposed versions of ProxSkip-VR algorithms on tomographic reconstruction problems (ProxSGDSkip, ProxSAGASkip, ProxSVRGSkip, ProxLSVRGSkip), and compare them against (i)  non-skipped  counterparts (ProxSGD, ProxSVRG, ProxLSVRG) (ii) stochastic methods without variance reduction (ProxSGDSkip), and (iii) deterministic full-gradient algorithms, skipped and non-skipped, (ISTA, FISTA, ProxSkip). We demonstrate a speed-up of at least $5\times$ of ProxSkip-VR algorithms over their non-skipped counterparts which so far have been considered among the most efficient methods with respect to computational time for  variational regularisation problems in tomography.
 
\section{ProxSkip-VR for Tomography}
\label{sec:proposed_method}
We consider two experimental setups to assess performance under different regularisations and evaluation criteria.

\subsection{TV tomography reconstruction on real dataset}\label{subsec:TVtomography}
Our goal is to compare instances of Algorithm \ref{alg:all_algos} in terms of the CPU time needed to reach a prescribed tolerance $\varepsilon>0$ to a high-accuracy reference solution $\bx^\ast$.
We solve \begin{equation}
\min_{\bx\in\mathbb{R}^n} \; \frac{1}{2}\|\bA\bx - \bb\|_{2}^{2}
+ \alpha \mathrm{TV}(\bx)
+ \mathbb{I}_{\{\bx \ge 0\}}(\bx),
\label{eq:tomography_problem}
\end{equation} 
where $\bA$ is the discrete Radon transform, on a real tomographic dataset, where ground truth is not available.  We use the isotropic TV as the regularisation term weighted with a regularisation parameter $\alpha>0$, optimised by visual inspection in the absence of ground truth. Here $\mathrm{prox}_{\alpha \mathrm{TV}}(\cdot)$, line 8 of Algorithm \ref{alg:all_algos}, is solved using the Fast Gradient Projection (FGP) algorithm\cite{BeckTeboulle2009} with a fixed number of iterations. 

We use a fixed number of warm-started FGP inner iterations to emulate cheap (10 iterations) and expensive (100 iterations) proximal steps, enabling controlled comparisons.\\[0.25em]
\textbf{Dataset:} We use a chemical imaging tomography dataset, see Figure \ref{fig:chemical_dataset_data_splitting}. The initial dataset was acquired for 800 projection angles with 695$\times$695 detector size and 700 vertical slices. For demonstration purposes we consider one 2D vertical slice with half the projections and $2\times$ rebinned detector size, i.e.\ 400 angles and 347 detector points. \\[0.5em]
\noindent\textbf{Reference solution and target accuracy:} The goal in this experiment is to identify which method achieves the desired reconstruction accuracy the fastest. The methods are evaluated in terms of the wall-clock CPU time required to either reach an error tolerance $\frac{\|\bx - \bx^*\|^{2}}{\|\bx^{*}\|^2}<10^{-5}$ or to reach 200 iterations for deterministic methods and 200 data passes for stochastic methods. Here, a data pass denotes one iteration for deterministic algorithms and when using $N$ subset gradients for stochastic algorithms. 
The high-accuracy reference solution $\bx^*$ is computed using a separate primal--dual proximal method to avoid bias, 
see Figure \ref{fig:chemical_dataset_data_splitting}.
\\[0.25em]
\noindent\textbf{Controlling gradient and proximal costs:} Our proposed algorithms allow independent control of the two dominant per-iteration costs.
(i) \emph{Proximal cost:} 
We  control the expected proximal workload through the skipping probability 
$p\in\{0.01, 0.05, 0.1, 0.3, 0.5\}$ which determines how often the (cheap or expensive) TV proximal operators are performed. (ii) \emph{Gradient cost:} For stochastic methods we vary the number of subsets $N$ in the decomposition of the data fidelity term. Larger $N$ corresponds to smaller-size subsets and therefore cheaper (but noisier) stochastic gradients. Decreasing $N$ yields more expensive gradients that better approximate the full gradient $\nabla f$ and have lower variance. 

\subsection{PnP-BM3D on a simulated dataset}\label{subsec:pnpbm3d}
We perform the tomographic reconstruction within a plug-and-play (PnP) framework using the BM3D denoiser and combining variance-reduced and skipping techniques \cite{PnP_VR}. \\[0.25em]
\noindent\textbf{Dataset:} We use simulated data, enabling quantitative evaluation of the iterates against a known ground truth. Specifically, we generate a cylindrical foam phantom with non-overlapping bubble microstructures \cite{Ching2017}, see Figure \ref{fig:BM3D_solutions_ssim_psnr}. The image dimensions match those of the real dataset in Section \ref{subsec:TVtomography}. We then compute analytical projection data from 400 uniformly spaced view angles and corrupt the measurements with additional simulated noise.\\[0.25em]
\noindent\textbf{BM3D denoiser:} 

We replace step 8 in Algorithms \ref{alg:all_algos} with a BM3D denoiser $D_{\sigma}(\cdot)$, with $\sigma>0$ separately optimised for all the algorithms.\\[0.25em]
\noindent\textbf{Stopping criterion and evaluation:} In this experiment, algorithms are compared with respect to the CPU time required to reach a given image-quality level, measured by PSNR and SSIM against the ground truth. We impose a fixed time budget and stop all methods after 3 minutes of CPU time, enabling a direct time-to-quality comparison.

\section{Numerical Results}
\begin{figure}
    \centering
    \includegraphics[width=0.98\linewidth]{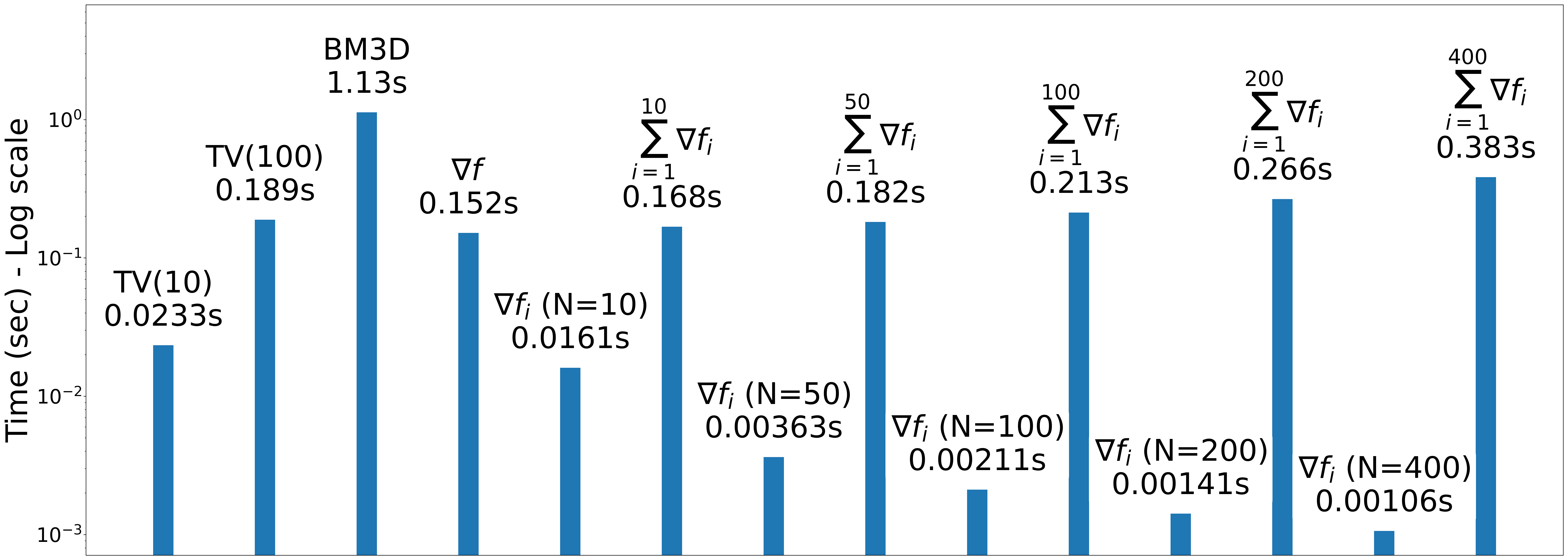}
    \caption{CPU time for TV proximal (10/100 iterations), BM3D, deterministic/stochastic gradients, and full stochastic gradients for different $N$.}
    \label{fig:time_comparison}
\end{figure}

All experiments were run on an Apple M2 Pro with 16 GB RAM and no GPU acceleration.
We first examine the average wall-clock time (10 runs) of the regularisers in each setting; see Figure~\ref{fig:time_comparison}. The TV-proximal is relatively inexpensive, taking 0.023sec and 0.189sec on average for 10 and 100 inner iterations respectively. In contrast, BM3D is substantially more costly, requiring more than 1sec per iteration.

\begin{figure}[h!]
    \centering
    \includegraphics[width=0.98\linewidth]{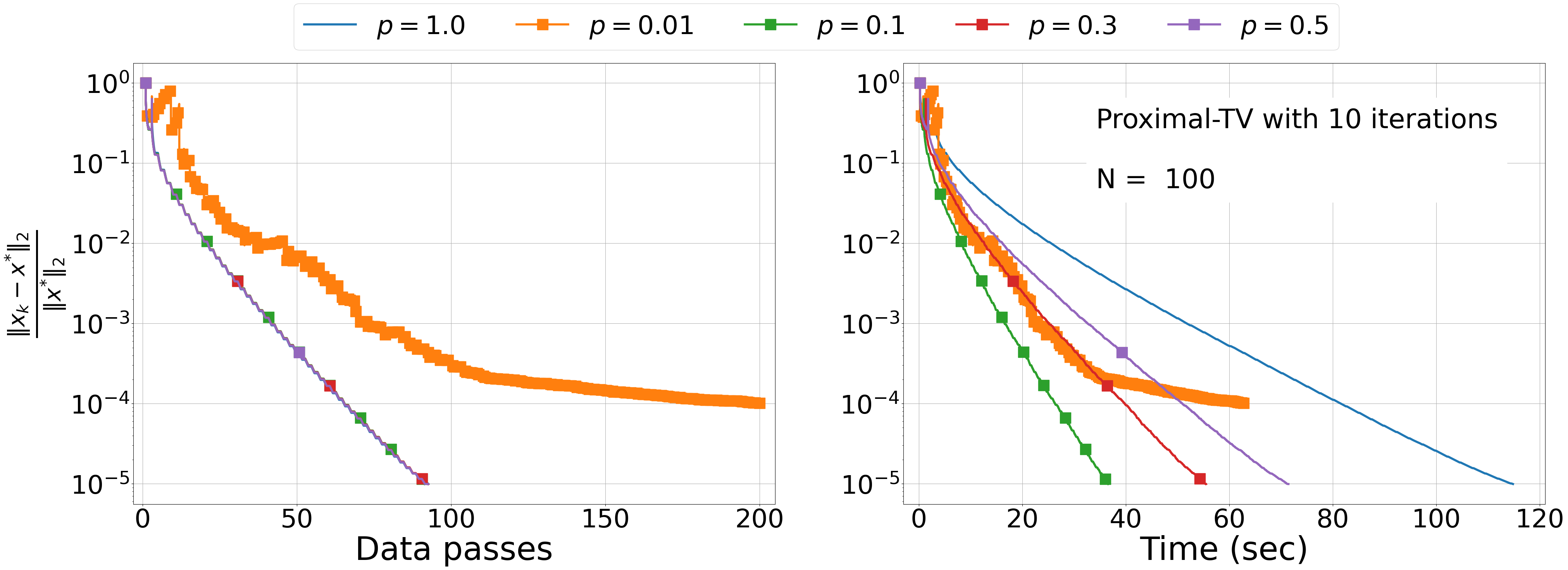}\\[0.5em]
    \includegraphics[width=0.98\linewidth]{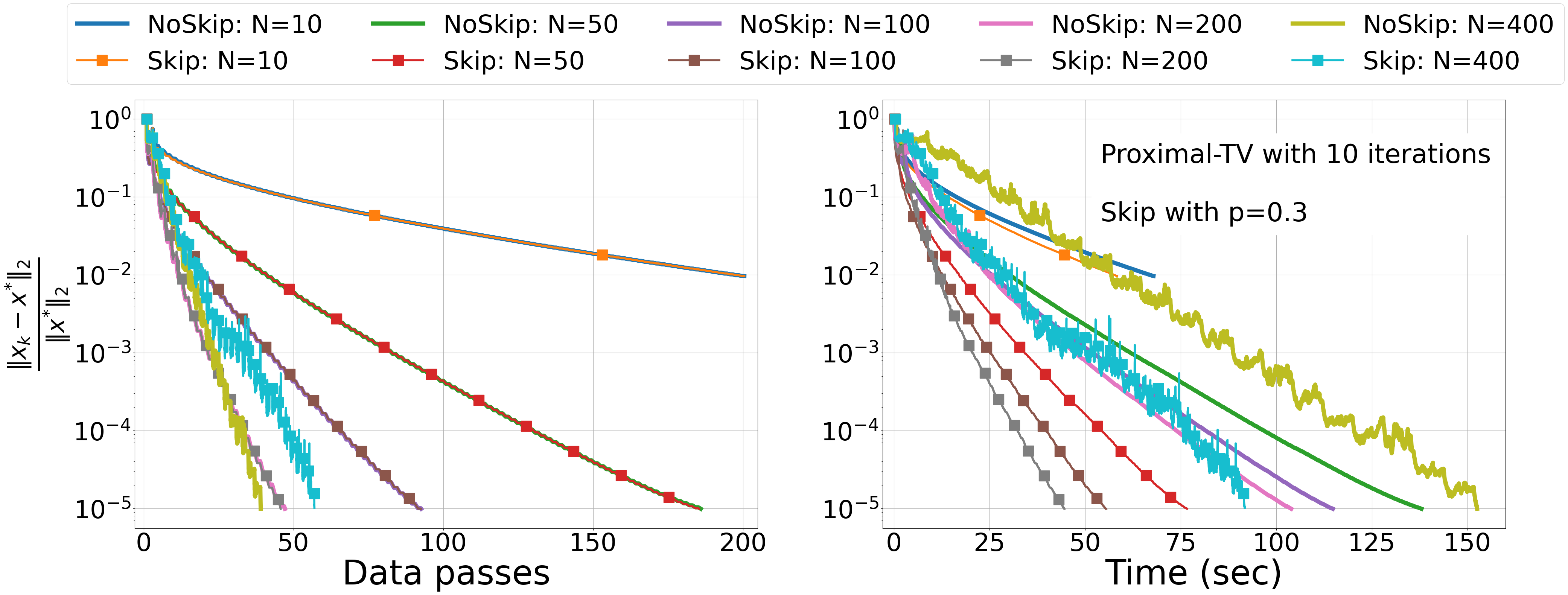}\\[0.5em]
    \vspace{0.6em}
    \begingroup
    \setlength{\tabcolsep}{3pt}
    \renewcommand{\arraystretch}{0.8}
    \centering
    \begin{tabular}{lccccc}
    \toprule
    ProxSVRG / $p$  & 0.01 & 0.1 & 0.3 & 0.5 & 1 \\
    \midrule
     & -- & \textbf{36.6\,sec} & 55.3\,sec & 71.1\,sec & 115\,sec \\
    \bottomrule
    \end{tabular} 
    
    \vspace{1em}    
    
    \begin{tabular}{lccccc}
    \toprule
    ProxSVRG / $N$ & 10 & 50 & 100 & 200 & 400 \\
    \midrule
    NoSkip & -- & 138\,sec & 115\,sec & 104\,sec & 153\,sec \\
    Skip   & -- & \textbf{76.6\,sec} & \textbf{55.3\,sec} & \textbf{44.6\,sec} & \textbf{91.8\,sec} \\
    \bottomrule
    \end{tabular}
    \endgroup
    \caption{ProxSVRG vs ProxSVRGSkip comparison in terms of data passes (left) and time (right). Top:   fixed number of subsets $N$ and varying probabilities $p$. Bottom: Fixed probability $p$ and
varying number of subsets N.}
    \label{fig:proxsvrg_vs_proxsvrg_skip_fix_N_fix_prob}
\end{figure}

We also report the cost of the deterministic full gradient $\nabla f$ used in every iteration for ISTA, FISTA and ProxSkip, the stochastic subset gradients $\nabla f_i$ used in the data-splitting cases, and the full-gradient computations required by ProxSVRG/ProxLSVRG for  different numbers of subsets $N$. For the stochastic gradients, we observe that increasing $N$ (using smaller subsets) substantially reduces the per-update cost of $\nabla f_i$. In contrast, the cost of a full gradient pass constructed from $N$ subset evaluations increases with $N$,  mainly due to the overhead of looping over more subsets. However, such full-gradient computations are performed infrequently in ProxSVRG (periodically) and randomly in ProxLSVRG, so they do not dominate runtime. Finally, the timings are consistent so that 
the cost of approximately $N$ stochastic gradients evaluations is roughly the same as full-gradient evaluation. \\[0.25em]
\noindent
\textbf{TV-Real dataset:} We first compare ProxSVRG and ProxSVRGSkip with  $N=100$ and $p\in\{0.01, 0.05, 0.1, 0.3, 0.5\}$,  Figure \ref{fig:proxsvrg_vs_proxsvrg_skip_fix_N_fix_prob} (top). The step size is fixed, $\gamma=\frac{1}{L}$. Similarly to observations in \cite{Papoutsellis2025}, the error trajectories in terms of the data passes are identical for all probabilities except the extreme case of $p=0.01$ in which the proximal is rarely applied. For all the other probabilities, we observe a faster convergence in terms of computational times, with a $3.13\times$ acceleration for $p=0.05$ compared to the non-skipped version and only 232 $\mathrm{prox}_{\alpha\mathrm{TV}}$ evaluations with 10 inner iterations.

In Figure \ref{fig:proxsvrg_vs_proxsvrg_skip_fix_N_fix_prob} (bottom), we perform the same comparison as above by keeping  the probability fixed, $p=0.3$  and varying the number of subsets $N$. When $N$ is small we observe that both skip and non-skip versions have slow convergence in terms of data passes and fail to reach high accuracy within 200 data passes. As $N$, and hence the variance, increase, both skip and non-skip versions converge in less than 200 data passes. Moreover, ProxSVRGSkip consistently reaches the target accuracy faster than ProxSVRG, with a $1.65\times$--$2.33\times$ speed-up in wall-clock time. The error trajectories plotted against data passes are identical for $N=50,100,200$ and begin to differ for $N=400$ due to the high stochastic variance. Interestingly, measured in data passes, ProxSVRG with $N=400$ reaches the target accuracy in the fewest data passes, but is the slowest in terms of computation time. This highlights that data passes do not fully capture the overall computational effort, particularly for large $N$, where each update is cheap but the stochastic gradient variance is high. 

\begin{table}[t]
\caption{Best time-to-accuracy for TV reconstruction (10 inner iterations):  Time ($T$) and number of iterations ($K$) at error thresholds $\varepsilon$; “--” denotes not reached.}
\label{tab:best_tv_10}
\hspace*{6mm} 
\centering
\begingroup
\small
\setlength{\tabcolsep}{3pt}
\renewcommand{\arraystretch}{0.8}

\begin{tabular}{lrrcc}
\toprule
Algorithm & $N$ & $p$ & $T/K~(\varepsilon=10^{-3})$ & $T/K~(\varepsilon=10^{-5})$ \\
\midrule
ISTA          & 1   & 1    & 182.1 / 1050 & 405.72 / 2336 \\
FISTA         & 1   & 1    & 72.9 / 415   & -- \\
ProxSkip      & 1  & 0.1  & 162.7 / 1037 & 373.6 / 2328 \\
ProxSAGA      & 100 & 1    & 128.9 / 6226 & -- \\
ProxSVRG      & 200 & 1    & 47.4 / 2113  & 103.9 / 4651 \\
ProxLSVRG     & 200 & 1    & 52.8 / 2419  & 109.8 / 4971 \\
ProxSAGASkip  & 200 & 0.05 & 14.9 / 6248  & -- \\
ProxSVRGSkip  & 200 & 0.05 & 11.3 / 2447  & \textbf{20.8 / 4545} \\
ProxLSVRGSkip & 200 & 0.05 & \textbf{11.1 / 2267}  & 21.7 / 4769 \\
\bottomrule
\end{tabular}
\endgroup
\end{table}

\begin{figure}[h!]
    \centering
    \begin{subfigure}[b]{0.48\textwidth}
        \centering
        \includegraphics[width=\textwidth]{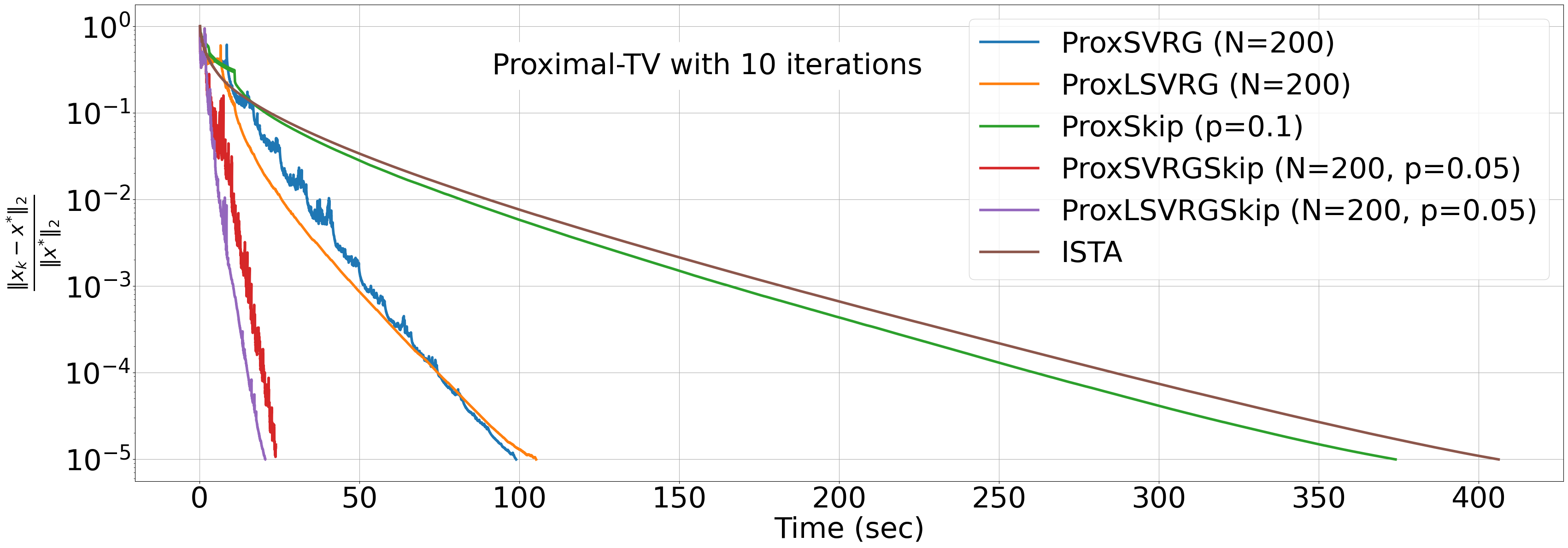}      
    \end{subfigure}
    \hfill
    \begin{subfigure}[b]{0.48\textwidth}
        \centering
        \includegraphics[width=\textwidth]{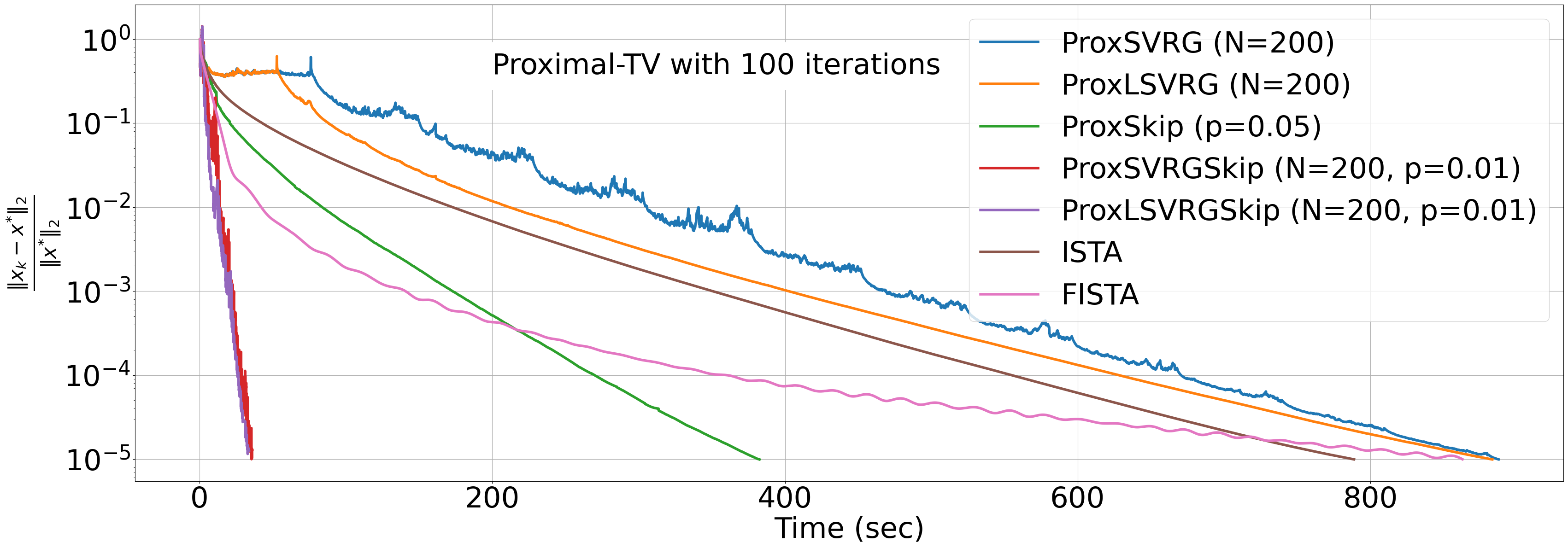}
    \end{subfigure}
    \caption{Convergence (relative error vs.\ CPU time) for TV tomography reconstruction with 10 (top) and 100 (bottom) inner iterations. Only methods reaching the target accuracy $10^{-5}$ within our computational budget are shown.}
    \label{fig:all_comparison_proxTV_10_100}
\end{figure}

We finally compare all the methods: ISTA ($\gamma=1.99/L$), FISTA ($\gamma=1/L$), stochastic estimators, ProxSAGA ($\gamma=1/3L$) and ProxSVRG/ProxLSVRG  ($\gamma=1/L$), and their skipped counterparts ProxSkip ($\gamma=1.99/L$), ProxSAGASkip ($\gamma=1/3L$), ProxSVRGSkip/ProxLSVRGSkip ($\gamma=1/L$). 
We evaluate all algorithms over the values of $N$ and  $p$ specified above (with full-gradient update probability set to $\frac{1}{N}$ where applicable) and for the two cases of inner iterations. In Figure \ref{fig:all_comparison_proxTV_10_100}, we report the best-performing  algorithms, showing only those that reached the target accuracy $\varepsilon=10^{-5}$ within our computational budget. We also report in Table~\ref{tab:best_tv_10} the best-performing configurations (10 inner iterations setting) at two target tolerances, $\varepsilon=10^{-3}$ and $10^{-5}$, summarising the fastest time-to-accuracy and the corresponding iteration count $K$ achieved by each method.

We observe \emph{striking} speed-ups.  ProxSVRGSkip and ProxLSVRGSkip achieve about a $5\times$ speed-up over their non-skipped counterparts and $20\times$ speed-up over deterministic ISTA. For 100 inner iterations, the corresponding speed-ups are about $23\times$ (vs. non-skipped) and $21\times$ (vs. ISTA). \\[0.25em]
\begin{figure}[h!]
    \centering
    \begin{subfigure}[b]{0.43\textwidth}
        \centering
        \includegraphics[width=\textwidth]{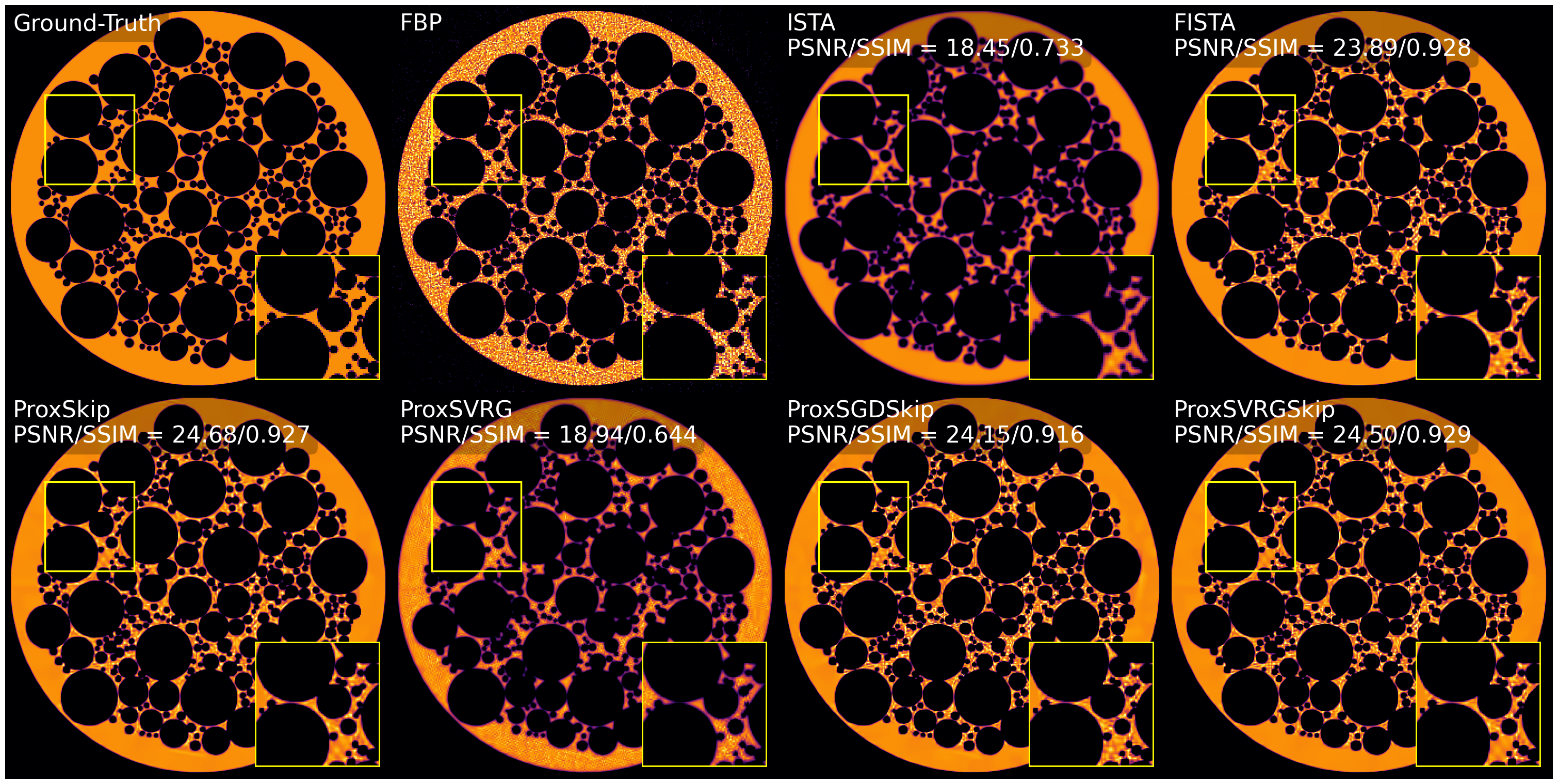}
    \end{subfigure}
    \hfill
    \begin{subfigure}[b]{0.48\textwidth}
        \centering
        \vspace{0.4em}
        \includegraphics[width=\textwidth]{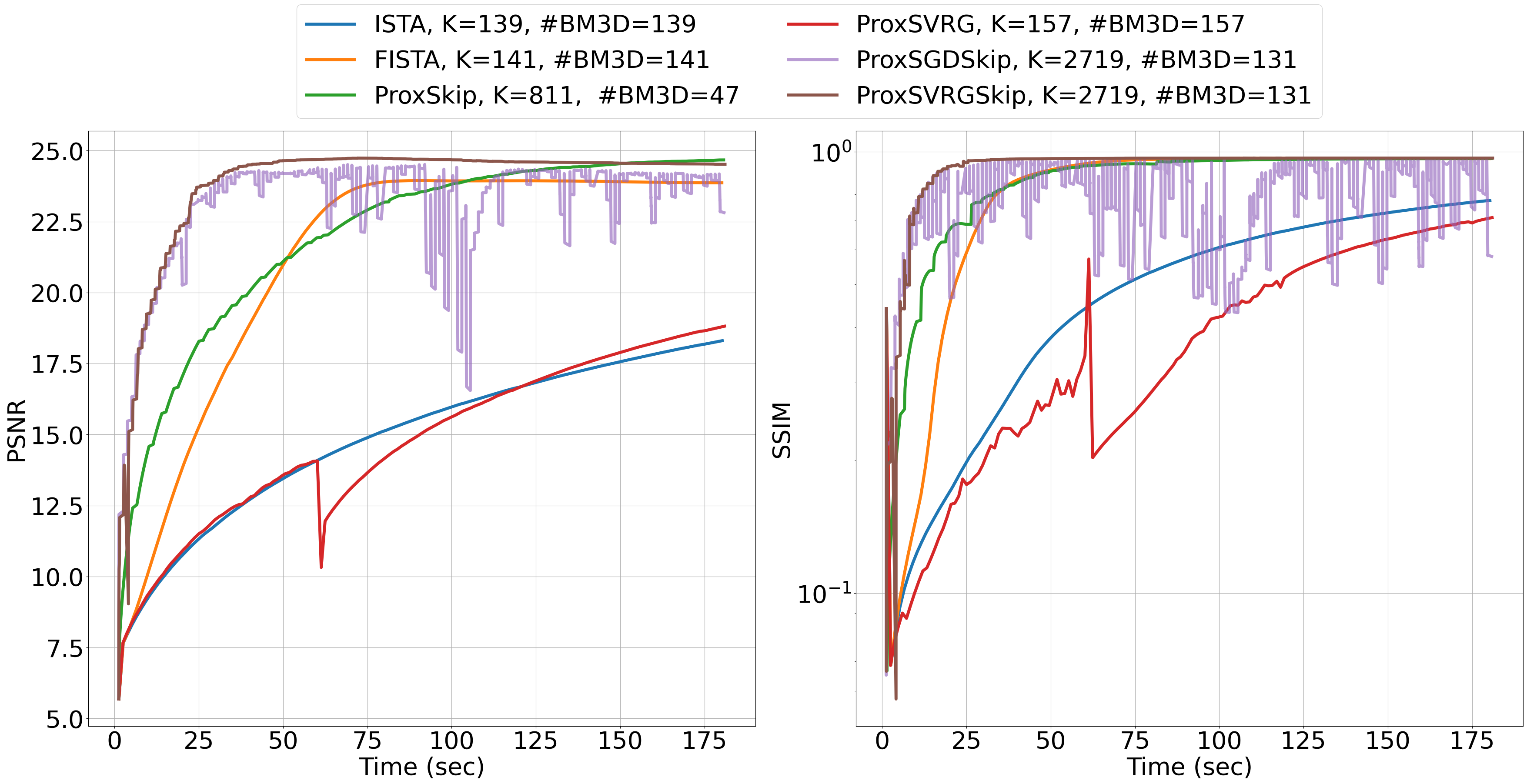}
    \end{subfigure}
    \caption{PnP–BM3D tomographic reconstruction on simulated foam data (3-minute budget). Top: ground truth, FBP, and reconstructions from 
    all
    algorithms. Middle: PSNR (left) and SSIM (right) versus wall-clock time. The number of subsets is $N=50$ with $p=0.05$ and $\#$BM3D denotes the number of calls for the BM3D denoiser.
    }
    \label{fig:BM3D_solutions_ssim_psnr}   
\end{figure}
\noindent
\textbf{PnP--BM3D-Simulated dataset: }In Figure~\ref{fig:BM3D_solutions_ssim_psnr}, we report the PnP--BM3D reconstructions obtained after a fixed CPU-time budget of 3 minutes, alongside the ground truth and the Filtered Back Projection (FBP) reconstruction. To ensure a fair comparison in this PnP setting, we tune the denoiser strength separately for skipped and non-skipped variants so that they converge to similar-quality reconstructions. For the non-skip algorithms, we use  $\sigma_{\text{non-skip}}=2\times10^{-4}$. For the skipped variants, we empirically found a simple heuristic that yields comparable image-quality metrics across all tested values of $p$, i.e.\ $\sigma_{\text{skip}} = \sigma_{\text{non-skip}}/\sqrt{p}$.
Since BM3D is approximately $42\times$ and $5\times$  more expensive than the cheap and expensive versions of  $\mathrm{prox}_{\mathrm{TV}}$, respectively, the overall per-iteration runtime is dominated by the denoiser rather than the gradient computation. As a result, and to keep the experimental sweep computationally feasible, we restrict the stochastic experiments to a single representative configuration $N=50, p=0.05$. 

The results show that skipping-based algorithms, e.g., ProxSkip and ProxSVRGSkip reach high PSNR/SSIM substantially faster than their non-skipped counterparts. 
For ISTA, FISTA, and ProxSVRG, BM3D is applied at every iteration, explaining their slower improvement in wall-clock time despite comparable iteration counts. Note, however, that ProxSVRG completes only three data passes for $N=100$ within the 3-minute budget, which explains its noticeably lower PSNR/SSIM in this experiment.\\[0.4em] 
\textbf{Code reproducibility:}
The code and datasets needed to reproduce the results will be available on the corresponding author's webpage.

\section{Conclusion}
Until now, variational regularisation problems in tomographic imaging have been most efficiently addressed using stochastic variance-reduced methods. Here we show for the first time that the performance of these algorithms can be significantly  accelerated when the proximal operator that corresponds to regularisation is randomly skipped. This combination of subset-based variance-reduced gradients with random proximal skipping, significantly reduces the overall computational cost while maintaining stable convergence. Our work opens up several directions for future research.  Integrating variance reduction and proximal skipping  to accelerated (momentum-based) schemes, has the potential to further reduce computational times. Furthermore, the incorporation of these mechanisms to stochastic primal-dual approaches is also worth of investigation. Finally, the adoption of ProxSkip-VR schemes to other challenging tasks, e.g.\ diffraction tomography and dynamic imaging, is also expected to be largely beneficial.

\newpage
\bibliographystyle{IEEEbib}
\bibliography{strings,refs}

\end{document}